
\tolerance=10000
\magnification=1200
\raggedbottom

\baselineskip=15pt
\parskip=1\jot

\def\sk{\vskip 3\jot}

\def\heading#1{\vskip3\jot{\noindent\bf #1}}
\def\label#1{{\noindent\it #1}}


\def\ref#1;#2;#3;#4;#5.{\item{[#1]} #2,#3,{\it #4},#5.}
\def\refinbook#1;#2;#3;#4;#5;#6.{\item{[#1]} #2, #3, #4, {\it #5},#6.} 
\def\refbook#1;#2;#3;#4.{\item{[#1]} #2,{\it #3},#4.}


\def\({\bigl(}
\def\){\bigr)}


\def\al{\alpha}
\def\be{\beta}
\def\ga{\gamma}

\def\ze{\zeta}
\def\et{\eta}

\def\io{\iota}
\def\ka{\kappa}

\def\ps{\psi}

\def\Ga{\Gamma}
\def\De{\Delta}

\def\Ph{\Phi}
\def\Ps{\Psi}



\def\Ex{{\rm Ex}}

\def\dbrack#1#2{\left[#1\brack#2\right]}

\def\Vol{{\rm Vol}}

{
\pageno=0
\nopagenumbers
\rightline{\tt cyclation.tex}
\vskip1in

\centerline{\bf Random Cyclations}
\vskip0.5in

\centerline{Nicholas Pippenger}
\centerline{\tt njp@princeton.edu}
\vskip0.25in

\centerline{Department of Computer Science}
\centerline{Princeton University}
\centerline{35 Olden Street}
\centerline{Princeton, NJ 08540}
\vskip1in

\noindent{\bf Abstract:}
Consider $n$ unit intervals,
say $[1,2], [3,4], \ldots,  [2n-1,2n]$.
Identify their endpoints in pairs at random, with all
$(2n-1)!! = (2n-1)(2n-3)\cdots 3\cdot 1$ pairings being equally likely.
The result is a collection of cycles of various lengths,
and we investigate the distribution of these lengths.
The distribution is similar to that of the distribution of the lengths
of cycles in a random permutation, but it also exhibits some striking 
differences.
\vfill\eject
}

\heading{1.  Introduction}

Our concern in this paper is with a process for randomly joining intervals into cycles.
Consider $n$ unit intervals,
say $[1,2], [3,4], \ldots,  [2n-1,2n]$.
Identify their endpoints in pairs at random, with all
$(2n-1)!! = (2n-1)(2n-3)\cdots 3\cdot 1$ pairings being equally likely.
The result, which we call a  {\it random $n$-cyclation}, is a collection of cycles of various lengths,
and we investigate the distribution of these lengths.

This process may be considered as a model for a chemical process of cyclation,
whereby $n$ molecules of a bilaterally symmetric monomer such as ethane, ${\rm C_2 H_6}$, are dehydogenated to form
cyclic polymers such as cyclic polyethylene, ${\rm (CH_2)}^{2n}$.
(The assumption that all pairings of endpoints are equally likely is not realistic, since
the likelihood of endpoints to pair will be influenced by their physical proximity.
We are, however, following a tradition in combinatorics going back to Cayley [C], in 
using an idealized model
to give our calculations a chemical interpretation.
The paper just cited, for example, overcounts the isomers of the paraffin series, since the trees enumerated
may have an exponentially growing number of vertices within a given distance of a particular vertex, whereas
a molecule in three-dimensional space may have only a cubically growing number of atoms within a given
distance of a particular atom.)

If the intervals to be joined were oriented, with for example $2m-1$ being the tail and $2m$ being the head
of the interval $[2m-1, 2m]$, and if they had to be joined head-to-tail consistently with this orientation, there
would then be just $n! = n \cdot (n-1)\cdots 2\cdot 1$ consistent pairings, and the distribution we seek
would be that of the lengths of cycles in a {\it random $n$-permutation}.
(See Goncharov [G1, G2] for pioneering work on this distribution.)
The distribution of cycle lengths for random cyclations is similar to that 
for random permutations, but it also exhibits some striking 
differences.

In Section 2 we shall study random cyclations primarily as regards their {\it number\/} of cycles,
largely disregarding the {\it lengths\/} of these cycles.
(In this respect, we are studying the {\it topology\/} of random cyclations, rather than their {\it geometry}.)

The number of $n$-permutations with $k$ cycles is a Stirling number of the first kind, which
(following Knuth [K2]) we denote ${n \brack k}$.
If
$$\Ph_n(\xi) = \sum_{1\le k\le n} {n\brack k} \, \xi^k$$
denotes the generating function for $n$-permutations according to their number of cycles,
then we have
$$\Ph_n(\xi) = \xi\,(\xi+1)\cdots(\xi+n-1) \eqno(1.1)$$
(see Knuth [K2], 1.2.9, (27)).
Setting $\xi=1$ yields the obvious
$$\sum_{1\le k\le n} {n\brack k} = n!.$$

Let us denote by $\dbrack{n}{k}$ the number of $n$-cyclations with $k$ cycles, and by
$$\Ps_n(\et) = \sum_{1\le k\le n} \dbrack{n}{k} \, \et^k$$
the genrating function for $n$-cyclations according to their number of cycles.
We shall show in Section 2 that
$$\eqalignno{
\Ps_n(\et) &= 2^n \, \Ph_n(\et/2) &(1.2) \cr
&= \et\,(\et+2)\cdots(\et+2n-2). &(1.3) \cr}$$
We shall first give a combinatorial proof of (1.2) by establishing the identity
$$\dbrack{n}{k} = 2^{n-k} \, {n\brack k}. \eqno(1.4)$$
We shall also present an analytic argument (inspired by a result of Harer and Zagier [H]) that gives
(1.3) directly, and also yields the following curious analogy: $n$-cyclations are to $n$-permutations as
the $(n-1)$-sphere is to the $(n-1)$-simplex.
In any case, setting $\et=1$ in (1.3) yields the obvious
$$\sum_{1\le k\le n} \dbrack{n}{k} = (2n-1)!!.$$

These results allow us to calculate the probabilities of the extreme cycle structures.
For the largest possible number $n$ of cycles, we have ${n\brack n} = \dbrack{n}{n} = 1$ from (1.1) and (1.3), 
so the probability that a random $n$-permutation or $n$-cyclation has $n$ cycles is $1/n!$ or $1/(2n-1)!!$,
respectively.
For the smallest possible number $1$ of cycles, we encounter our first surprise: from (1.1) and (1.3)
we have ${n\brack 1} = (n-1)!$ and $\dbrack{n}{1} = 2^{n-1} \, (n-1)!$.
Thus the probability that a random $n$-permutation or $n$-cyclation has just a single cycle is $1/n$ or
$2^{n-1} \, (n-1)!/(2n-1)!!$, respectively.
We have
$$\eqalign{
{2^{n-1} \, (n-1)!\over (2n-1)!!}
&= {2^{2n-1} \, (n-1)! \, n!\over (2n)!} \cr
&\sim {1\over 2}\,\left({\pi \over n}\right)^{1/2} \cr}$$
(where we have used the identity $(2n-1)!! = (2n)!/2^n\,n!$ and Stirling's asymptotic formula).
Thus the probability of having just a single cycle is substantially larger for cyclations than for permutations.

These results also allow us to calculate the expected number of cycles (as well as the higher moments, if desired)
in random permutations and cyclations.
From (1.1) and (1.3), we have
$$\eqalign{
{\Ph'_n(1)\over n!}
&= \sum_{1\le m\le n} {1\over m}\cr
&= H_n \cr
&= \log n + \ga + O\left({1\over n}\right), \cr
}$$
where $\ga = 0.577\ldots$ is Euler's constant (see Knuth [K2], 	1.2.7, (3)), and
$$\eqalign{
{\Ps'_n(1)\over (2n-1)!!}
&= \sum_{1\le m\le n} {1\over 2m-1}\cr
&= H_{2n} - \textstyle{1\over 2}H_n \cr
&= {\textstyle}{1\over 2}\log n + \log 2  + {\ga\over 2} + \displaystyle O\left({1\over n}\right). \cr
}$$
(All logarithms in this paper are to the natural base $e = 2.718\ldots\,$.)

In Section 3 we shall study the lengths of the cycles in random cyclations.
Our approach is modelled upon that used by Shepp and Lloyd [S] to study the cycle lengths in random permutations.
We summarize their approach as follows.

From an exact expression (formula (2.1) below)
for the number of $n$-permutations with cycles structure
$(i_1, i_2, \ldots, i_n)$ (that is, with $i_1$ cycles of length $1$, $i_2$ cycles of length $2$, $\ldots\,$,
and $i_n$ cycles of length $n$), it is apparent that the cycle structure of a random $n$-permutation has
the following interpretation:
let $i_1, i_2, \ldots, i_n$ be independent Poisson-distributed random variables with means
$1, 1/2, \ldots, 1/n$, respectively, and condition on the event  that
cycles with these lengths contain altogether $n$ elements (that is, that
$n = \sum_{1\le l\le n} l\, i_l$).
Equivalently, we may 
let $i_1, i_2, \ldots, i_l, \ldots$ be independent Poisson-distributed random variables with means
$1, 1/2, \ldots, 1/l, \ldots\,$, respectively, and condition on the event that $n = \sum_{l\ge 1} l\, i_l$).
(This event now has probability $0$, since it implies that $i_{n+1} = i_{n+2} = \cdots = 0$.)
Shepp and Lloyd hit upon the idea of letting 
$i_1, i_2, \ldots, i_l, \ldots$ be independent Poisson-distributed random variables,
without any conditioning.
If this were done with the means indicated above, the sum $\nu =  \sum_{l\ge 1} l\, i_l$, representing the
total number of elements in all the cycles, would be infinite with probability $1$.
Thus they fix $0<z<1$ and use the means $z, z^2/2, \ldots, z^l/l, \ldots\,$.
This results in $\nu$ having a geometric distribution: 
the probability that $\nu = n$ is $(1-z)\,z^n$, and $(1-z)/(1-z\ze)$ is the generating function for $\nu$.
Furthermore, when $\nu = n$, the conditional distribution of the cycle structure
$(i_1, i_2, \ldots, i_n)$ is that of a random $n$-permutation.
(Since $n=0$ is possible in this model, we should agree that there is just one $0$-permutation,
which of course has no cycles.)
This allows asymptotic information about the cycle structure as $n\to\infty$ to be obtained
by applying Tauberian theorems to results concerning the limit $z\to 1$.
Shepp and Lloyd show, among other things, that the expectation of the length $L_n$ of the longest cycle in a random
$n$-permutation has the asymptotic behavior
$$\Ex[L_n] \sim \int_0^\infty e^{-E(x)-x} \,dx \,n,$$
where
$$E(x) = \int_x^\infty {e^{-y}\over y} \, dy,$$
and
$\int_0^\infty e^{-E(x)-x} \,dx = 0.6243\ldots\,$,
while the expectation of the length $S_n$ of the shortest cycle in a random $n$-permutation satisfies
$$\Ex[S_n] \sim e^{-\ga} \, \log n,$$
where $e^{-\ga} = 0.5614\ldots\,$.
(We note that the function $E(x)$ can be expressed as the incomplete gamma-function $\Ga(0,x)$,
or  in terms of the exponential integral, as $-{\rm Ei}(-x)$; see Lebedev [L], (3.1.3).)

For $n$-cyclations, we have a similar interpretation, with the means of the 
independent Poisson-distributed random variables being $1/2, 1/4, \ldots, 1/2l, \ldots$ before conditioning.
We then let the means be $z/2, z^2/4, \ldots, z^l/2l, \ldots$ without conditioning.
The number $\nu = \sum_{l\ge 1} l\, i_l$ of intervals then has a negative-binomial distribution:
the probability that $\nu = n$ is $(1-z)^{1/2} \, {2n\choose n}\,(z/4)^n$, 
and $\((1-z)/(1-z\ze)\)^{1/2}$ is the generating function for $\nu$. 
Furthermore, when $\nu = n$, the conditional distribution of the cycle structure
$(i_1, i_2, \ldots, i_n)$ is that of a random $n$-cyclation.
(We agree of course that there is just one $0$-cyclation, and thus that $(-1)!! = 1$.)
Again applying Tauberian theorems, we encounter our second surprise.
We shall show in Section 3 that the expectation of the length $M_n$ of the longest cycle in a random
$n$-cyclation satisfies
$$\Ex[M_n] \sim \int_0^\infty e^{-E(x)/2-x} \,dx \,n, \eqno(1.5)$$
where
$$\int_0^\infty e^{-E(x)/2-x} \,dx = 0.7578\ldots,$$
while the expectation of the length $T_n$ of the shortest cycle in a random $n$-cyclation satisfies
$$\Ex[T_n] \sim {\pi^{1/2}\over 2}\int_0^\infty e^{E(x)/2-x} \,dx \,n^{1/2}, \eqno(1.6)$$
where $\pi = 3.1415\ldots$ is the circular ratio, and
$(\pi^{1/2}/2)\int_0^\infty e^{E(x)/2-x} \,dx = 1.4572\ldots\,$.
Thus the expected length of the shortest cycle is substantially larger for cyclations than for permutations.
\sk

\heading{2. The Number of Cycles}

In this section we shall give two proofs of (1.3).
Each proof in fact establishes the connection (1.2) between cyclations and permutations.

The first proof esablishes the identity (1.4).
We begin by deriving a formula for the number ${n\brack k}$ of $n$-permutations with $k$ cycles.
Consider the $n$-permutations with cycle structure $(i_1, i_2, \ldots, i_n)$,
where $n = \sum_{1\le l\le n} l\,i_l$ and $k = \sum_{1\le l\le n} i_l$.
Such a permutation can be specified as follows.
Arrange the elements $1, 2, \ldots, n$ in some order from left to right (which can be done in $n!$ ways).
Partition these elements into blocks by forming $i_1$ blocks of size $1$ from the leftmost $i_1$ elements,
$i_2$ blocks of size two from the next $2\,i_2$ elements, and so forth, with the last $n\,i_n$ elements
forming $i_n$ blocks of size $n$.
Form a cycle of the permutation from each block by letting the successor of each element of the block
except the rightmost be the next element to the right, with the successor of the rightmost element 
of the block being the leftmost element of that block.
It is clear that every $n$-permutation with cycle structure  $(i_1, i_2, \ldots, i_n)$ can be obtained in this way.
Indeed, since the same $n$-permutation is obtained if the 
block of elements forming a cycle of length $l$ are cyclically
permuted among themselves
(which can be done in $l$ ways), or if the blocks of elements forming the cycles of length $l$ are 
permuted among themselves (which can be done in $i_l!$ ways), there are exactly
$$\prod_{1\le l\le n} l^{i_l} \, i_l!$$
left-to-right orders that yield any given $n$-permutation with cycle structure $(i_1, i_2, \ldots, i_n)$.
Thus if $[i_1, i_2, \ldots, i_n]$ denotes the set of  $n$-permutations with cycle structure
$(i_1, i_2, \ldots, i_n)$ (it is in fact a conjugacy class in the symmetric group on $n$ elements),
then we have
$$\#[i_1, i_2, \ldots, i_n] = {n!\over \displaystyle\prod_{1\le l\le n} l^{i_l} \, i_l!}, \eqno(2.1)$$
and thus
$${n\brack k} = 
\sum_{(i_1, i_2, \ldots, i_n)} {n!\over \displaystyle\prod_{1\le l\le n} l^{i_l} \, i_l!}, \eqno(2.2)$$ 
where the sum is over all cycle structures $(i_1, i_2, \ldots, i_n)$ such that $n = \sum_{1\le l\le n}
l\,i_l$ and $k = \sum_{1\le l\le n} i_l$.

Next we derive an analogous formula for the number $\dbrack{n}{k}$ of $n$-cyclations with $k$ cycles.
Consider the $n$-cyclations with cycle structure $(i_1, i_2, \ldots, i_n)$,
where $n = \sum_{1\le l\le n} l\,i_l$ and $k = \sum_{1\le l\le n} i_l$.
Such a cyclation can be specified as follows.
Arrange the elements $1, 2, \ldots, n$ in some order from left to right (which can be done in $n!$ ways),
and add primes to some subset of these elements (which can be done in $2^n$ ways).
Replace each unprimed element $m$ by the interval $[2m-1,2m]$ with its smaller endpoint $2m-1$ to the left
of its larger endpoint $2m$, and 
replace each primed element $m'$ by the interval $[2m-1,2m]$ with its larger endpoint to the left
of its smaller endpoint.
Partition these intervals into blocks by forming $i_1$ blocks of size $1$ from the leftmost $i_1$ intervals,
$i_2$ blocks of size two from the next $2\,i_2$ intervals, and so forth, with the last $n\,i_n$ intervals
forming $i_n$ blocks of size $n$.
Form a cycle of the cyclation from each block by letting the  right endpoint
of each interval of the block
except the rightmost be paired with the left endpoint of the 
next interval to the right, with the right endpoint  of the rightmost interval 
of the block being paired with the left endpoint of the leftmost interval of that block.
It is clear that every $n$-cyclation with cycle structure  $(i_1, i_2, \ldots, i_n)$ can be obtained in this way.
Indeed, since the same $n$-cyclation is obtained if the 
elements in a block giving rise to a cycle of length $l$ are cyclically
permuted among themselves
(which can be done in $l$ ways), or if they are 
cyclically permuted after reversing their left-to-right order and complementing their primed-or-unprimed status
(which can be done in another $l$ ways),
or if the blocks of elements giving rise to  the cycles of length $l$ are 
permuted among themselves (which can be done in $i_l!$ ways), there are exactly
$$\prod_{1\le l\le n} (2l)^{i_l} \, i_l!$$
left-to-right orders and choices of primed elements
that yield any given $n$-cyclation with cycle structure $(i_1, i_2, \ldots, i_n)$.
Thus if $[[i_1, i_2, \ldots, i_n]]$ denotes the set of  $n$-cyclations with cycle structure
$(i_1, i_2, \ldots, i_n)$,
then we have
$$\#[[i_1, i_2, \ldots, i_n]] = {n!\, 2^n\over \displaystyle\prod_{1\le l\le n} (2l)^{i_l} \, i_l!}, \eqno(2.3)$$
and thus
$$\dbrack{n}{k} = 
\sum_{(i_1, i_2, \ldots, i_n)} {n!\,2^n\over \displaystyle\prod_{1\le l\le n} (2l)^{i_l} \, i_l!}, \eqno(2.4)$$ 
where the sum is over all cycle structures $(i_1, i_2, \ldots, i_n)$ such that $n = \sum_{1\le l\le n}
l\,i_l$ and $k = \sum_{1\le l\le n} i_l$.

Comparing (2.2) with (2.4) yields (1.4).
It is easy to see that this proof of (1.4) is ``bijective'', in the sense that the derivations of (2.2) and (2.4)
can be used to give an explicit combinatorial correspondence between appropriately defined sets of size
$2^n\,{n\brack k}$ and $2^k\,\dbrack{n}{k}$.
Multiplying by $\et^k$ and summing over $1\le k\le n$ then yields (1.2).
Finally, from (1.2) and the definition of $\Ph_n(\xi)$ we obtain (1.3).

We shall now give a second, analytic rather than bijective, proof of (1.3).
The starting point for this proof is a rather circuitous proof of  (1.1).
We begin by observing that since both sides of (1.1) are polynomials in $\xi$, it will suffice to establish
it for $\xi$ a positive integer.
To this end, we shall consider the number of $n$-permutations in which each cycle is colored with one of $\xi$ 
colors.
Since a permutation with $k$ cycles can be colored in $\xi^k$ ways, the number of such colored $n$-permutations
is $\Ph_n(\xi)$.
We can count the number of colored permutations in another way as follows.
We begin by assigning a color to each of the elements $1, 2, \ldots, n$.
We then arrange the elements into cycles in such a way that all the elements in a given cycle have the same color.
We then sum the number of permutations obtained over all ways of assigning colors to elements.
There are ${n\choose j_1,j_2,\ldots,j_\xi}$ ways of assigning the first color to $j_1$ elements, the second
color to $j_2$ elements, and so forth, with $j_\xi$ elements receiving the $\xi$-th color.
(We must of course have $n = \sum_{1\le\al\le\xi} j_\al$.)
Given such a distribution of colors, there are $j_\al!$ ways of arranging the $j_\al$ elements receiving color $\al$
into cycles (that is, permuting them among themselves).
Thus
$$\Ph_n(\xi) = 
\sum_{(j_1, j_2, \ldots, j_\xi)} {n\choose j_1,j_2,\ldots,j_\xi} \; \prod_{1\le\al\le\xi} j_\al!, \eqno(2.5)$$
where the sum is over all compositions $(j_1, j_2, \ldots, j_\xi)$ of $n$ into $\xi$ non-negative parts.
From the formula
$${n\choose j_1,j_2,\ldots,j_\xi} = {n!\over \displaystyle\prod_{1\le\al\le\xi} j_\al!}, \eqno(2.6)$$
we obtain
$$\eqalign{
\Ph_n(\xi) 
&= \sum_{(j_1, j_2, \ldots, j_\xi)} n! \cr
&= {n+\xi-1\choose \xi-1} \; n! \cr
&= \xi\,(\xi+1)\cdots(\xi+n-1),}$$
which is (1.1), since the number of compositions of $n$ into $\xi$ non-negative parts is 
${n+\xi-1\choose \xi-1} = (n+\xi-1)!/(\xi-1)!\,n!$.

If we did not know the formula (2.6), we could still proceed as follows.
Replace each factorial in (2.5) with its integral representation 
$j! = \int_0^\infty x^j \, e^{-x} \, dx$.
This yields
$$\eqalign{
\Ph_n(\xi)
&=\sum_{(j_1, j_2, \ldots, j_\xi)} {n\choose j_1,j_2,\ldots,j_\xi} \cr
&\qquad\qquad\int_0^\infty \int_0^\infty \cdots \int_0^\infty 
x_1^{j_1} x_2^{j_2} \cdots x_\xi^{j_\xi} \,
e^{-(x_1 + x_2 + \cdots + x_\xi)} \, dx_1 \, dx_2 \cdots dx_\xi \cr
&= \int_0^\infty \int_0^\infty \cdots \int_0^\infty 
(x_1 + x_2 + \cdots + x_\xi)^n \,
e^{-(x_1 + x_2 + \cdots + x_\xi)} \, dx_1 \, dx_2 \cdots dx_\xi, \cr
}$$
by the multinomial theorem.
Since the integrand is now a function of $s = x_1 + x_2 + \cdots + x_\xi$, we obtain
$$\eqalign{
\Ph_n(\xi) 
&= \int_0^\infty s^n \, e^{-s} \;\; s^{\xi-1}  \, ds \;\; \Vol_{\xi-1}(\De_{\xi-1}) \cr
&= (n+\xi-1)! \, \Vol_{\xi-1}(\De_{\xi-1}) \cr
&= \xi\,(\xi+1)\cdots(\xi+n-1), \cr
}$$
where 
$$\Vol_{\xi-1}(\De_{\xi-1}) = 1/(\xi-1)! \eqno(2.7)$$
denotes the $(\xi-1)$-dimensional volume of the 
$(\xi-1)$-dimensional simplex $\De_{\xi-1}$ defined by the equation $x_1 + x_2 + \cdots + x_\xi = 1$ and 
the inequalities $x_1\ge 0, x_2\ge 0, \ldots, x_\xi\ge 0$ in $\xi$-dimensional space.
Thus we again obtain (1.1).
(We observe that the identity (2.7) can be derived in the following way.
Multiplying $\xi$ copies of the identity
$1 = \int_0^\infty e^{-x} \, dx$ we obtain
$$1 = \int_0^\infty \int_0^\infty \cdots  \int_0^\infty e^{-(x_1 + x_2 + \cdots + x_\xi)} \, 
dx_1 \, dx_2\cdots dx_\xi.$$ 
This multiple integral can be evaluated like the preceding one,
$$\eqalign{
\int_0^\infty \int_0^\infty \cdots  \int_0^\infty e^{-(x_1 + x_2 + \cdots + x_\xi)} \, 
dx_1 \, x_2 \cdots dx_\xi
&= \int_0^\infty e^{-s} \;\;  s^{\xi-1}\, ds \;\;
\Vol_{\xi-1}(\De_{\xi-1}) \cr
&= (\xi-1)! \, \Vol_{\xi-1}(\De_{\xi-1}). \cr
}$$
Combining these two evaluations of the multiple integral yields (2.7).)

Next we give an analogous derivation of (1.3).
It will suffice to establish it for $\et$ an even positive integer.
By considering the number of ways of assigning one of $\et$ colors to each cycle of each 
$n$-cyclation we obtain
$$\Ps_n(\et) = \sum_{(j_1, j_2, \ldots, j_\et)} {n\choose j_1,j_2,\ldots,j_\et} \; 
\prod_{1\le\be\le\et} (2j_\be - 1)!!.$$
We now use the integral representation 
$(2j-1)!! = (1/2\pi)^{1/2}\,\int_{-\infty}^\infty y^{2j} \, e^{-{1\over 2}y^2}\,dy$ to obtain
$$\eqalign{
\Ps_n(\et)
&={1\over(2\pi)^{\et/2}}\sum_{(j_1, j_2, \ldots, j_\et)} {n\choose j_1,j_2,\ldots,j_\et} \cr
&\qquad\qquad\qquad\int_{-\infty}^\infty \int_{-\infty}^\infty \cdots \int_{-\infty}^\infty 
y_1^{2j_1} y_2^{2j_2} \cdots y_\et^{2j_\et} \,
e^{-{1\over 2}(y^2_1 + y^2_2 + \cdots + y^2_\et)} \, dy_1 \, dy_2 \cdots dy_\et \cr
&= {1\over(2\pi)^{\et/2}}\int_{-\infty}^\infty \int_{-\infty}^\infty \cdots \int_{-\infty}^\infty 
(y^2_1 + y^2_2 + \cdots + y^2_\et)^n \,
e^{-{1\over 2}(y^2_1 + y^2_2 + \cdots + y^2_\et)} \, dy_1 \, dy_2 \cdots dy_\et. \cr
}$$
Since the integrand is now a function of $r = (y^2_1 + y^2_2 + \cdots + y^2_\et)^{1/2}$, we obtain
$$\eqalign{
\Ps_n(\et) 
&= {1\over(2\pi)^{\et/2}}\int_0^\infty r^{n/2} \, e^{-{1\over 2}r^2} \;\; r^{\et-1} \,  dr \;\;
\Vol_{\et-1}(S_{\et-1}) \cr
&= {2^{n-1}\over\pi^{\et/2}}\int_0^\infty s^n \, e^{-s} \;\; s^{\et/2-1} \,  ds \;\;
\Vol_{\et-1}(S_{\et-1}) \cr 
&= {2^{n-1} \,  ({\et/2} + n - 1)!\over\pi^{\et/2}} \;\;
\Vol_{\et-1}(S_{\et-1}) \cr 
&= {2^n \, ({\et/2} + n - 1)!\over (\et/2 - 1)!} \cr
&= \et\,(\et+2)\cdots(\et+2n-2), \cr
}$$
where we have made the substitution $r = (2s)^{1/2}$, and
$$\Vol_{\et-1}(S_{\et-1}) = {2\,\pi^{\et/2} \over  (\et/2 - 1)!} \eqno(2.8)$$
denotes the $(\et-1)$-dimensional volume of the 
$(\et-1)$-dimensional sphere $S_{\et-1}$ defined by the equation $y^2_1 + y^2_2 + \cdots + y^2_\et = 1$  
in $\et$-dimensional space.
Thus we again obtain (1.3).
(We observe that the identity (2.8) can be derived in the following way.
Multiplying $\et$ copies of the identity
$(2\pi)^{1/2} = \int_{-\infty}^\infty e^{-{1\over 2}y^2} \, dy$ we obtain
$$(2\pi)^{\et/2} = \int_{-\infty}^\infty \int_{-\infty}^\infty \cdots  \int_{-\infty}^\infty 
e^{-{1\over 2}(y^2_1 + y^2_2 + \cdots + y^2_\et)} \, 
dy_1 \, dy_2\cdots dy_\et.$$ 
This multiple integral can be evaluated like the preceding one,
$$\eqalign{
\int_{-\infty}^\infty \int_{-\infty}^\infty \cdots  \int_{-\infty}^\infty 
e^{-{1\over 2}(y^2_1 + y^2_2 + \cdots + y^2_\et)} \, 
&dy_1 \, dy_2\cdots dy_\et \cr
&= \int_0^\infty e^{-{1\over 2}r^2} \;\;  r^{\et-1}\, dr \;\;
\Vol_{\et-1}(S_{\et-1}) \cr
&= 2^{\et/2-1} \,\int_0^\infty e^{-s} \;\;  s^{\et/2-1}\, ds \;\;
\Vol_{\et-1}(S_{\et-1}) \cr
&= 2^{\et/2-1} \, (\et/2-1)! \, \Vol_{\et-1}(S_{\et-1}). \cr
}$$
Combining these two evaluations of the multiple integral yields (2.8).)
\sk

\heading{3. The Lengths of Cycles}

In this section we shall derive (1.5) and (1.6), establishing the asymptotic behaviors of the expectations 
$\Ex[M_n]$ and $\Ex[T_n]$ of the lengths
$M_n$ and $T_n$
of the longest and shortest cycles, in a random  $n$-cyclation as $n\to\infty$.
This will be done in three steps.
First, we shall introduce a new notion of random cyclation, we we call a 
random $z$-cyclation, where $0<z<1$.
For a random $z$-cyclation, the total number of intervals is itself a random variable $\nu$.
If however we condition on the event $\nu=n$, the conditional distribution is that of a random 
$n$-cyclation.
Second, we establish the asymptotic behaviors of the expectations $\Ex[M_z]$ and $\Ex[T_z]$
of the lengths
$M_z$ and $T_z$
of the longest and shortest cycles, in a random $z$-cyclation as $z\to 1$.
Finally, we apply Tauberian theorems to deduce from these results the asymptotic behaviors of 
$\Ex[M_n]$ and $\Ex[T_n]$ as $n\to \infty$.

For the first step, let $z$ be a real number in the range $0<z<1$.
We shall define a {\it random $z$-cyclation\/} as follows.
First, let $\io_1, \io_2, \ldots, \io_l, \ldots$ be independent Poisson-distributed random variables with
means $z/2, z^2/4, \ldots, z^l/2l, \ldots\,$, respectively.
Next, define the random variable $\nu$ by
$$\nu = \sum_{l\ge 1} l\,\io_l.$$
Finally, choose a $\nu$-cyclation with cycle structure $(\io_1, \io_2, \ldots, \io_\nu)$, with all
$$\#[[\io_1, \io_2, \ldots, \io_\nu]] = 
{\nu!\, 2^\nu \over \displaystyle\prod_{1\le l\le \nu} (2l)^{\io_l} \, \io_l!} \eqno(3.1)$$
cyclations having this cycle structure being equally likely.
(Since $n$ will always denote an integer, while $z$ never will, no confusion between the notions
of random $n$-cyclations and random $z$-cyclations should be possible.)

The number $\ka = \sum_{l\ge 1} \io_l$ of cycles in a random $z$-cyclation is a sum of 
independent Poisson-distributed random variables, and thus is itself Poisson-distributed, with mean
$$\eqalign{
\Ex[\ka] 
&= \sum_{l\ge 1} {z^l\over 2l} \cr
&= {1\over 2} \log {1\over 1-z}. \cr}$$

Of more importance for us will be the distribution of the random variable $\nu$.
For any cycle structure $(i_1, i_2, \ldots, i_l, \ldots)$
(where we agree to identify cycle structures that differ from one another by appending or deleting
any number, finite or infinite, of $0$s at the right-hand end), we have
$$\Pr[(\io_1, \io_2, \ldots, \io_l, \ldots) = (i_1, i_2, \ldots, i_l, \ldots)] = 
{(1-z)^{1/2} \, z^n \over \displaystyle\prod_{1\le l\le n} (2l)^{i_l} \, i_l!}, \eqno(3.2)$$
where
$$n = \sum_{l\ge 1} l\,i_l. \eqno(3.3)$$
Thus we have
$$\Pr[\nu = n] = (1-z)^{1/2} \, z^n \; \sum_{(i_1, i_2, \ldots, i_n)} 
{1 \over \displaystyle\prod_{1\le l\le n} (2l)^{i_l} \, i_l!}, \eqno(3.4)$$
where the sum is over all cycle structures $(i_1, i_2, \ldots, i_l, \ldots)$ satisfying (3.1).
To evaluate the sum in (3.4), we sum (2.3) over all cycle structures satisfying (3.1).
The result is
$$(2n-1)!! = n! \, 2^n \; \sum_{(i_1, i_2, \ldots, i_n)} 
{1 \over \displaystyle\prod_{1\le l\le n} (2l)^{i_l} \, i_l!}.$$
Combining this with (3.4) yields
$$\eqalignno{
\Pr[\nu = n]
&= (1-z)^{1/2} \, z^n \; {(2n-1)!! \over n! \, 2^n} \cr
&= (1-z)^{1/2} \, z^n \; {1\over 4^n}{2n\choose n}. &(3.5)\cr
}$$
Thus $\nu$ has a negative-binomial distribution: the generating function for $\nu$ is
$$\ps_z(\ze) = \left({1-z\over 1-z\ze}\right).$$

Finally, we consider the conditional distribution of a random $z$-cyclation, given the event $\nu = n$.
By comparing the probability (3.2) of a given cycle structure with the number (3.1) of equally likely
cyclations with that cycle structure, we see that all $(2n-1)!!$ $n$-cyclations have equal probability
$(1-z)^{1/2} \, z^n /  n! \, 2^n$.
Thus, the distribution of a random $z$-cyclation, conditioned on the event $\nu = n$, is that
of a random $n$-cyclation.

In particular, from (2.5) we obtain
$$\eqalignno{
\Ex[M_z]
&= \sum_{n\ge 0} \Pr[\nu=n] \, \Ex[M_n] \cr
&= (1-z)^{1/2} \; \sum_{n\ge 0} {z^n\over 4^n} {2n\choose n} \, \Ex[M_n]. &(3.6)\cr
}$$

The next step is to determine the asymptotic behavior of $\Ex[M_z]$ as $z\to 1$.
We begin by determining the probability of the event $M_z = l$.
For $l\ge 1$, define
$$t_l = \sum_{1\le m\le l-1} {z^m \over 2m},$$
with the understanding that $t_1 = 0$ and
$$\eqalign{
t_\infty 
&= \sum_{1\le m\le \infty} {z^m \over 2m} \cr
&= {1\over 2} \log {1\over 1-z}. \cr}$$
The event $M_z = l$ occurs if and only if 
(1) $\io_l \ge 1$, which occurs with probability $1 - e^{-(t_{l+1} - t_l)}$, and
(2) $\io_{l+1} = \io_{l+2} = \cdots = 0$, which occurs independently
with probability $e^{-(t_\infty - t_{l+1})}$.
Thus we have
$$\eqalignno{
\Pr[M_z = l]
&= e^{-(t_\infty - t_{l+1})} - e^{-(t_\infty - t_l)} \cr
&= \int_{t_l}^{t_{l+1}} e^{t_\infty - t} \, dt. &(3.7)\cr}$$

From (3.7) we have
$$\eqalign{
\Ex[M_z] 
&= \sum_{l\ge 1} l \, \Pr[M_z = l] \cr
&= \sum_{l\ge 1} l \; \int_{t_l}^{t_{l+1}} e^{t_\infty - t} \, dt. \cr}$$
Making the substitution 
$$t_\infty - t = E(x)/2, \eqno(3.8)$$
we obtain
$$\Ex[M_z]  = {1\over 2}\sum_{l\ge 1} l \; \int_{x_l}^{x_{l+1}} {e^{-E(x)/2 - x}\over x}\,dx, \eqno(3.9)$$
where $x_l$ is the unique solution of the equation
$$t_\infty - t_l = E(x_l)/2.$$
(The existence and uniqueness of $x_l$,
together with the inequalities $0<x_1<x_2<\cdots<x_l<x_{l+1}<\cdots<x_\infty = \infty$,
follow from the fact that as $l$ increases from $1$ to $\infty$,
$t_\infty - t_l$ strictly decreases from $t_\infty$ to $0$, while as
$x$ increases from $0$ to $\infty$,
$E(x)/2$ strictly decreases from $\infty$ to $0$.)
Thus (3.9) may be rewritten
$$\Ex[M_z]  =  {1\over 2}\int_{x_1}^\infty e^{-E(x)/2 - x}\;{l(x)\over x}\,dx, \eqno(3.10)$$
where $l(x)$ denotes the unique integer $l$ such that $x_l \le x < x_{l+1}$.

To determine the asymptotic behavior of the integral in (3.10), we shall need an estimate for $l(x)$.
For $y$ in the range $-m\log z < y < -(m+1)\log z$, we have
$${z^{m+1}\over -(m+1)\log z} < {e^{-y}\over y} < {z^m \over -m\log z}.$$
Integrating over $y$ in an integral of length $-\log z$ yields
$${z^m\over m+1} < \int_{-m\log z}^{-(m+1)\log z} {e^{-y}\over y} \, dy < {z^m \over m},$$
and summing over $m$ yields
$$E(-l\log z) < t_\infty - t_l < E\(-(l+1)\log z\),$$
which is equivalent to 
$$-(l-1)\log z < x_l < -l\log z.$$
Since $-\log z = (1-z) + O\((1-z)^2\)$ as $z\to 1$,
we have
$${l(x)\over x} = {1\over 1-z} + O(1) \eqno(3.11)$$
in this limit.

We shall also need an estimate for $x_1$.
From the formula
$$E(x) =  - \log x - \ga + O(x)$$
as $x\to 0$
(see Lebedev [L], (3.1.5)),
and the definition
$$\eqalign{
{E(x_1) / 2}
&= t_\infty - t_1 \cr
&= {1\over 2} \log {1\over 1-z}, \cr}$$
we obtain 
$$x_1 = e^{-\ga} \, (1-z) + O\((1-z)^2\) \eqno(3.12)$$
as $z\to 1$.

Applying the estimates (3.11) and (3.12) to (3.10), we have
$$\eqalignno{
\Ex[M_z]  
&=  {1\over 2} \int_{x_1}^\infty e^{-E(x)/2 - x}\,dx \, {1\over 1-z} + 
O\left( \int_{x_1}^\infty {e^{-E(x)/2 - x}\over x}\,dx \right) \cr
&=  {1\over 2} \int_{0}^\infty e^{-E(x)/2 - x}\,dx \, {1\over 1-z} + 
O\left( \log {1\over 1-z}\right), &(3.13) \cr
}$$
since
$$\eqalign{
\int_{x_1}^\infty e^{-E(x)/2 - x}\,dx
&= \int_{0}^\infty e^{-E(x)/2 - x}\,dx - \int_{0}^{x_1} e^{-E(x)/2 - x}\,dx \cr
&= \int_{0}^\infty e^{-E(x)/2 - x}\,dx + O(x_1) \cr
&= \int_{0}^\infty e^{-E(x)/2 - x}\,dx  + O(1-z) \cr
}$$
(the integrand being bounded by $1$), and
$$\eqalign{
 \int_{x_1}^\infty {e^{-E(x)/2 - x}\over x}\,dx
&= O\left(\log {1\over x_1}\right)\cr
&= O\left(\log {1\over 1-z}\right)\cr
}$$
(the integrand being bounded by $e^{-x}$ for $x\ge 1$ and by $1/x$ for $x\le 1$).
Formula (3.13) gives the desired estimate for $\Ex[M_z]$.

The final step is to deduce the asymptotic behavior of $\Ex[M_n]$ as $n\to\infty$ from that of
$\Ex[M_z]$ as $z\to 1$.
This will be done through the use of Tauberian theorems, which require ``side conditions'' for their application.
In the case at hand, these side conditions require that the coefficents 
${1\over 4^n} {2n\choose n} \, \Ex[M_n]$ in (3.6) be non-negative, which is trivial, and non-decreasing in $n$,
which is equivalent to
$$(2n+1)\,\Ex[M_{n+1}] \ge (2n+2)\,\Ex[M_n]. \eqno(3.14)$$
To verify this condition, we start with a random $n$-cyclation $C$ and randomly insert an $(n+1)$-st interval $I$,
obtaining a random $(n+1)$-cyclation.
This cannot decrease the length of the longest cycle, and if the insertion puts $I$ into a longest cycle of $C$
(which happens with probability at least $2\,\Ex[M_n]/(2n+1)$), the length of the longest cycle increases by $1$.
Thus we have
$$\Ex[M_{n+1}] \ge \Ex[M_n]  + {2\over 2n+1} \, \Ex[M_n],$$
or equivalently
$$(2n+1)\,\Ex[M_{n+1}] \ge (2n+3)\,\Ex[M_n],$$
which is even stronger that (3.14).

From (3.6) and (3.13) we have
$$\eqalign{
\sum_{n\ge 0} {z^n\over 4^n} {2n\choose n} \, \Ex[M_n]
&= {\Ex[M_z]\over (1-z)^{1/2}} \cr
&\sim {1\over 2}\int_{0}^\infty e^{-E(x)/2 - x}\,dx \, {1\over (1-z)^{3/2}} \cr
}$$
as $z\to 1$.
Since the coefficients ${1\over 4^n} {2n\choose n} \, \Ex[M_n]$ are non-negative,
the final Tauberian theorem in Section 7.5 of de Bruijn [B] can be applied with ($\ga = 3/2$) to give
$$\sum_{0\le n\le m} {1\over 4^n} {2n\choose n} \, \Ex[M_n] \sim 
{2\over 3\pi^{1/2}}  \int_{0}^\infty e^{-E(x)/2 - x}\,dx \, m^{3/2}$$
as $m\to\infty$.
Since the coefficients ${1\over 4^n} {2n\choose n} \, \Ex[M_n]$ are non-decreasing in $n$,
the summatory analog of the Tauberian theorem in Section 7.3 of de Bruijn [B] can be applied (with $\al = 1/2$)
to give
$${1\over 4^n} {2n\choose n} \, \Ex[M_n] \sim 
{2\over \pi^{1/2}} \int_{0}^\infty e^{-E(x)/2 - x}\,dx \, n^{1/2}$$
as $n\to \infty$.
(This summatory analog is in fact a corollary of the theorem in Section 7.3 of de Bruijn [B]:
given a non-negative and non-decreasing sequence $f_n$, define $f(t)$ for $t\ge 0$ by $f(t) = f_n$ for
$n\le t < n+1$; then $f(t)$ is also non-negative and non-decreasing;  if 
$F_m = \sum_{0\le n\le m}  \sim (\al+1)^{-1}\, m^{\al+1}$, then we also have 
$F(t) = \int_0^t f(x)\,dx \sim (\al+1)^{-1}\, t^{\al+1}$; the theorem cited then implies that
$f(t) \sim t^\al$, and thus we also have $f_n \sim n^\al$.)
Finally, since 
$${1\over 4^n} {2n\choose n} \sim {1\over \pi^{1/2} \, n^{1/2}},$$
we obtain (1.5) as desired.

We now turn to the expectation $\Ex[T_n]$ of the length $T_n$ of the shortest cycle in a random $n$-cyclation.
Using the notion of a random $z$-cyclation as before, we have the analog of (3.6):
$$\eqalignno{
\Ex[T_z]
&= \sum_{n\ge 0} \Pr[\nu=n] \, \Ex[T_n] \cr
&= (1-z)^{1/2} \; \sum_{n\ge 0} {z^n\over 4^n} {2n\choose n} \, \Ex[T_n]. &(3.15)\cr
}$$

Our next step is to determine the asymptotic behavior of $\Ex[T_z]$ as $z\to 1$.
We begin by determining the probability of the event $T_z = l$.
The event $T_z = l$ occurs if and only if 
(1) $\io_l \ge 1$, which occurs with probability $1 - e^{-(t_{l+1} - t_l)}$, and
(2) $\io_{1} = \io_{2} = \cdots \io_{l-1} = 0$, which occurs independently
with probability $e^{-t_l}$.
and therefore
$$\eqalignno{
\Pr[T_z = l]
&= e^{-t_l} - e^{-t_{l+1}} \cr
&= \int_{t_l}^{t_{l+1}} e^{ - t} \, dt, \cr}$$
This yields
$$\eqalign{
\Ex[T_z] 
&= \sum_{l\ge 1} l \, \Pr[T_z = l] \cr
&= \sum_{l\ge 1} l \; \int_{t_l}^{t_{l+1}} e^{ - t} \, dt. \cr}$$
Again making the substitution (3.8), and proceding as before,
we obtain
$$\eqalign{
\Ex[T_z]  
&= {(1-z)^{1/2}\over 2}\sum_{l\ge 1} l \; \int_{x_l}^{x_{l+1}} {e^{E(x)/2 - x}\over x}\,dx \cr
&=  {(1-z)^{1/2}\over 2}\int_{x_1}^\infty e^{E(x)/2 - x}\;{l(x)\over x}\,dx \cr
}$$
Applying the estimates (3.11) and (3.12), we have
$$\eqalignno{
\Ex[T_z]  
&=  {1\over 2} \int_{x_1}^\infty e^{E(x)/2 - x}\,dx \, {1\over (1-z)^{1/2}} + 
O\left((1-z)^{1/2} \int_{x_1}^\infty {e^{E(x)/2 - x}\over x}\,dx \right) \cr
&=  {1\over 2} \int_{0}^\infty e^{-E(x)/2 - x}\,dx \, {1\over (1-z)^{1/2}} + 
O\left((1-z)^{1/2} \log {1\over 1-z}\right). &(3.16) \cr
}$$
Formula (3.16) gives the desired estimate for $\Ex[T_z]$.

Our final step is to deduce the asymptotic behavior of $\Ex[T_n]$ as $n\to\infty$ from that of
$\Ex[T_z]$ as $z\to 1$.
This will again be done through the use of Tauberian theorems. 
In this case, the side conditions require that the coefficents 
${1\over 4^n} {2n\choose n} \, \Ex[T_n]$ in (3.15) be non-negative, which is trivial, and non-decreasing in $n$,
which is equivalent to
$$(2n+2)\,\Ex[T_n] \le (2n+1)\,\Ex[T_{n+1}]. \eqno(3.17)$$
To verify this condition, we start with a random $(n+1)$-cyclation $C$ and delete the $(n+1)$-st interval
$I$, obtaining a random $n$-cyclation. 
If $I$ is alone in a cycle of length $1$, which happens with probability $1/(2n+1)$,
the shortest cycle is that of a random $n$-cyclation.
With the remaining probability $2n/(2n+1)$, the shortest cycle is that of a random $(n+1)$-cyclation, except that
if $I$ is taken from a shortest cycle in $C$,
which happens with probability at least $2\Ex[T_{n+1}]/(2n+1)$, the length of the shortest cycle is decreased by $1$.
Thus we have
$$\Ex[T_n] \le {1\over 2n+1}\,\Ex[T_n] + {2n\over 2n+1}\,\Ex[T_{n+1}] - {2\over 2n+1}\,\Ex[T_{n+1}],$$
or equivalently
$$n\,\Ex[T_n] \le (n-1)\,\Ex[T_{n+1}],$$
which is even stronger than (3.17).

From (3.15) and (3.16) we have
$$\eqalign{
\sum_{n\ge 0} {z^n\over 4^n} {2n\choose n} \, \Ex[T_n]
&= {\Ex[T_z]\over (1-z)^{1/2}} \cr
&\sim {1\over 2} \int_{0}^\infty e^{E(x)/2 - x}\,dx \, {1\over 1-z} \cr
}$$
as $z\to 1$.
Since the coefficients ${1\over 4^n} {2n\choose n} \, \Ex[T_n]$ are non-negative,
the final Tauberian theorem in Section 7.5 of de Bruijn [B] can be applied with ($\ga = 1$) to give
$$\sum_{0\le n\le m} {1\over 4^n} {2n\choose n} \, \Ex[T_n] \sim 
{1\over 2}  \int_{0}^\infty e^{E(x)/2 - x}\,dx \, m$$
as $m\to\infty$.
Since the coefficients ${1\over 4^n} {2n\choose n} \, \Ex[T_n]$ are non-decreasing in $n$,
the summatory analog of the Tauberian theorem in Section 7.3 of de Bruijn [B] can be applied (with $\al = 0$)
to give
$${1\over 4^n} {2n\choose n} \, \Ex[T_n] \to 
{1\over 2} \int_{0}^\infty e^{E(x)/2 - x}\,dx$$
as $n\to \infty$.
Finally, since 
$${1\over 4^n} {2n\choose n} \sim {1\over \pi^{1/2} \, n^{1/2}},$$
we obtain (1.6) as desired.
\sk

\heading{4. Conclusion}

We have determined the distribution of the number of cycles, and the expectations of the lengths of the
longest and shortest cycles, in a random $n$-cyclation.
There are many directions in which further results might be found.
As was done for random permutations by Shepp and Lloyd [S], one might seek the 
expectations of the lengths of  the $r$-th longest and $r$-shortest cycles
(these are presumably asymptotic to appropriate constants times $n$ and $n^{1/2}(\log n)^{r-1}$,
respectively), or more generally  the $m$-th moments of these quantities.
One might also try to determine the actual distributions of the length of the  shortest cycle
(this presumably tends to the limiting distribution $\Pr[T_n = l]\to e^{-H_l/2} - e^{-H_{l+1}/2}$ as
$n\to\infty$ with $l$ fixed, where $e^{-H_l/2} - e^{-H_{l+1}/2} \sim e^{-\ga/2}/2l^{3/2}$ as $l\to\infty$),
or of the normalized length $M_n/n$ of the longest cycle, which presumably tends to that of 
the largest component in a Poisson-Dirichlet distribution with parameter $1/2$.
(Poisson-Dirichlet distributions were introduced by Kingman [K1].
The Poisson-Dirichlet distribution with parameter $1/2$ has been shown by Aldous [A] to govern the 
normalized size of the largest component in the graph of a random mapping from a set to itself.)
\sk

\heading{5.  References}
 
\refinbook A; D. J. Aldous;
``Exchangeability and Related Topics'';
in: P.~L. Hennequin (Ed.);
\'Ecole d'\'Et\'e de Probabilit\'es de Saint-Flour XIII---1983;
Lecture Notes in Mathematics, v.~1117, Springer-Verlag, 1985.

\refbook B; N. G. de Bruijn;
Asymptotic Methods in Analysis;
Dover Publications, 1981.

\ref C; A. Cayley;
``On the Mathematical Theory of Isomers'';
Phil.\ Mag.; 47 (1874) 444.

\ref G1; V. Goncharov;
``Sur la distribution des cycles dans les permutations'';
C. R. (Doklady) Acad.\ Sci.\ URSS (N. S.);
35 (1942) 267--269.

\ref G2; V. Goncharov;
``Du domaine d'analyse combinatoire'';
Bull.\ Acad.\ Sci.\ USSR Ser.\ Math.\ (Izv.\ Akad.\ Nauk SSSR);
8 (1944) 3--48;
{\it Amer.\ Math.\ Soc.\ Transl.\ (2)}, 19 (1962) 1--46.

\ref H; J. Harer and D. Zagier;
``The Euler Characteristic of the Moduli Space of Curves'';
Invent.\ Math.; 85 (1986) 457--485.

\ref K1; J. F. C. Kingman;
``Random Discrete Distributions'';
J. Roy.\ Stat.\ Soc.\ (B); 37 (1975) 1--22.

\refbook K2; D. E. Knuth;
The Art of Computer Programming, Volume 1: Fundamental Algorithms;
Addison-Wesley Publishing Company, 1968.

\refbook L; N. N. Lebedev;
Special Functions and Their Applications;
Prentice-Hall, Incorporated, 1965.

\ref S; L. A. Shepp and S. P. Lloyd;
``Ordered Cycle Lengths in a Random Permutation'';
Trans.\ Amer.\ Math.\ Soc.; 121 (1966) 340--357.

\bye